\documentclass[a4paper,11pt]{amsart}
\usepackage{graphicx}
 \usepackage{mathptmx}
\usepackage{amsmath}
\usepackage{amssymb}
\usepackage{enumitem}
\usepackage{xcolor}

\newmuskip\pFqmuskip

\newcommand*\pFq[6][8]{%
  \begingroup 
  \pFqmuskip=#1mu\relax
  \mathcode`=\string"8000
  \begingroup\lccode`\~=`\,
  \lowercase{\endgroup\let~}\pFqcomma
  F^{#2}_{#3}{\left(\genfrac..{0pt}{}{#4}{#5}\bigg|#6\right)}%
  \endgroup
}
\newcommand{\pFqcomma}{\mskip\pFqmuskip}

\newtheorem{theorem}{Theorem}[section]

\begin{document}

\title[Probabilistic degenerate Dowling polynomials associated with random variables]{Probabilistic degenerate Dowling polynomials associated with random variables}

\author{Taekyun  Kim}
\address{Department of Mathematics, Kwangwoon University, Seoul 139-701, Republic of Korea}
\email{tkkim@kw.ac.kr}
\author{Dae San  Kim }
\address{Department of Mathematics, Sogang University, Seoul 121-742, Republic of Korea}
\email{dskim@sogang.ac.kr}


\begin{abstract}
The aim of this paper is to study probabilistic versions of the degenerate Whitney numbers of the second kind and those of the degenerate Dowling polynomials, namely the probabilistic degenerate Whitney numbers of the second kind associated with $Y$ and the probabilistic degenerate Dowling polynomials associated with $Y$. Here Y is a random variable whose moment generating function exists in some neighborhood of the origin. We derive some properties, explicit expressions, certain identities, recurrence relations and generating functions for those numbers and polynomials. In addition, we investigate their generalizations, namely  the probabilistic degenerate $r$-Whitney numbers of the second kind associated with $Y$ and the probabilistic degenerate $r$-Dowling polynomials associated with $Y$, and get similar results to the aforementioned numbers and polynomials.
\end{abstract}

\maketitle

\section{Introduction}
Carlitz initiated a study of degenerate versions of special polynomials and numbers in his work on the degenerate Bernoulli and degenerate Euler polynomials and numbers (see [7]). It is remarkable that various degenerate versions of many special numbers and polynomials have been explored recently by some mathematicians not only with their number-theoretic or combinatorial interests but also with their applications to other areas, including probability, quantum mechanics and differential equations. In the course of this quest, many different tools are employed, like generating functions, combinatorial methods, $p$-adic analysis, umbral calculus, operator theory, differential equations, special functions, probability theory and analytic number theory (see [8,10,12-15,17,18] and the references therein). \par
Assume that $Y$ is a random variable such that the moment generating function of $Y$
\begin{equation}
E[e^{tY}]=\sum_{n=0}^{\infty}E[Y^{n}]\frac{t^{n}}{n!},\,\,(|t|<r) \,\,\,\mathrm{exists\,\, for\,\, some} \,\,r>0, \label{0}
\end{equation}
where $E$ stands for the mathematical expectation (see [3,4,16,24]).
The aim of this paper is to study probabilistic versions of the degenerate Whitney numbers of the second kind and those of the degenerate Dowling polynomials, namely the probabilistic degenerate Whitney numbers of the second kind associated with $Y$, $W_{m,\lambda}^{Y}(n,k)$ and the probabilistic degenerate Dowling polynomials associated with $Y$, $D_{m,\lambda}^{Y}(n,x)$. In addition, we investigate their generalizations, namely  the probabilistic degenerate $r$-Whitney numbers of the second kind associated with $Y$, $W_{m,\lambda}^{(Y,r)}(n,k)$ and the probabilistic degenerate $r$-Dowling polynomials associated with $Y$, $D_{m,\lambda}^{(Y,r)}(n,x)$. Then we derive some properties, explicit expressions, certain identities, recurrence relations and generating functions for those numbers and polynomials. \par
The outline of this paper is as follows. In Section 1, we recall the Whitney numbers of the second kind, the $r$-Whitney numbers of the second kind, the Dowling polynomials and the $r$-Dowling polynomials. We remind the reader of the degenerate exponentials and the degenerate Stirling numbers of the second. We recall the degenerate Whitney numbers of the second kind, the degenerate $r$-Whitney numbers of the second kind, the degenerate Dowling polynomials and the degenerate $r$-Dowling polynomials. We remind the reader of the partial Bell polynomials and the complete Bell polynomials. Then, for any random variable $Y$ satisfying the moment condition (see \eqref{0}), we recall the definition of the probabilistic degenerate Stirling numbers of the second kind associated with $Y$. Section 2 is the main result of this paper. Assume that $Y$ is a random variable satisfying the moment condition in \eqref{0}. Let $(Y_{j})_{j\ge 1}$ be a sequence of mutually independent copies of the random variable $Y$, and let $S_{k}=Y_{1}+\cdots+Y_{k},\ (k\ge 1)$, with $S_{0}=0$, \,\,(see\ [1,3,4,16,25,26]). We define $W_{m,\lambda}^{Y}(n,k)$. Then three different expressions of those numbers are obtained in terms of expectations of various random variables in Theorems 2.1-2.3. Another expression is derived in terms of the partial Bell polynomials in Theorem 2.15. In Theorem 2.4, we derive a finite sum identity involving those numbers and $E\Big[(mS_{k}+1)_{n,\lambda}\Big]$. Then we define $D_{m,\lambda}^{Y}(n,x)$. We obtain the generating function of $D_{m,\lambda}^{Y}(n,x)$ in Theorem 2.5 and two explicit expressions of those polynomials in Theorems 2.6 and 2.7. We derive a recurrence relation for $D_{m,\lambda}^{Y}(n,x)$ in Theorem 2.8. In Theorem 2.9, we get an identity for $D_{m,\lambda}^{Y}(n,x)$, which shows that those polynomials do not satisfy the binomial identity. We deduce a finite sum identity involving $D_{m,\lambda}^{Y}(n,x)$ and the partial Bell polynomials in Theorem 2.11. Moreover, higher order derivatives of $D_{m,\lambda}^{Y}(n,x)$ are obtained in Theorem 2.14. Next, we define $W_{m,\lambda}^{(Y,r)}(n,k)$. We get an explicit expression for those polynomials in Theorem 2.10. In Theorems 2.12 and 2.13, we derive identities involving those polynomials and the partial Bell polynomials. Then we define $D_{m,\lambda}^{(Y,r)}(n,x)$. Finally, we obtain the generating function and an explicit expression, respectively in Theorem 2.16 and Theorem 2.17. In the rest of this section, we recall the facts that are needed throughout this paper. \par

\vspace{0.1in}

A finite lattice $L$ is geometric if it is a finite semimodular lattice which is also atomic. Dowling constructed an important finite geometric lattice $Q_{n}(G)$ out of a finite set of $n$ elements and a finite group $G$ of order $m$, called Dowling lattice of rank $n$ over a finite group of order $m$ (see [11,13,15]).
If $L$ is the Dowling lattice $Q_{n}(G)$ of rank $n$ over a finite group of $G$ order $m$, then the Whitney numbers of the second kind $W_{Q_{n}(G)}(n,k)$ are denoted by $W_{m}(n,k)$ and satisfy the following relation
\begin{equation}
x^{n}	=\sum_{k=0}^{n}W_{m}(n,k)m^{k}\bigg(\frac{x-1}{m}\bigg)_{k},\quad (n\ge 0),\quad (\mathrm{see}\ [11,13,15]),\label{1}
\end{equation}
 where $(x)_{0}=1,\ (x)_{n}=x(x-1)\cdots(x-n+1),\ (n\ge 1)$. \par
 From \eqref{1}, we have
 \begin{equation}
 (mx+1)^{n}=\sum_{k=0}^{n}W_{m}(n,k)m^{k}(x)_{k},\quad (n\ge 0),\quad (\mathrm{see}\ [11,13,15,18]).\label{2}
 \end{equation}
Let $r$ be a nonnegative integer. Then, as generalizations of $W_{m}(n,k)$, the $r$-Whitney numbers of the second kind are defined by
\begin{equation}
(mx+r)^{n}=\sum_{k=0}^{n}W_{m}^{(r)}(n,k)m^{k}(x)_{k},\quad (n\ge 0),\quad (\mathrm{see}\ [13]).\label{3}
\end{equation} \par
For $n\ge 0$, the Dowling polynomials $D_{m}(n,x)$ and the $r$-Dowling polynomials $D_{m}^{(r)}(n,x)$ are respectively given by
\begin{equation}
D_{m}(n,x)=\sum_{k=0}^{n}W_{m}(n,k)x^{k},\quad (n\ge 0),\label{4}	
\end{equation}
and
\begin{displaymath}
D_{m}^{(r)}(n,x)=\sum_{k=0}^{n}W_{m}^{(r)}(n,k)x^{k},\quad (n\ge 0),\quad (\mathrm{see}\ [15]).
\end{displaymath} \par
It is well known that the Stirling numbers of the second kind ${n \brace k}$ are defined by
\begin{equation}
x^{n}=\sum_{k=0}^{n}{n \brace k}(x)_{k},\quad (n\ge 0),\quad (\mathrm{see}\ [1-25]). \label{5}
\end{equation}
For any nonzero $\lambda\in\mathbb{R}$, the degenerate exponentials are defined by
\begin{equation}
e_{\lambda}^{x}(t)=(1+\lambda t)^{\frac{x}{\lambda}}=\sum_{k=0}^{\infty}(x)_{k,\lambda}\frac{t^{k}}{k!},\quad (\mathrm{see}\ [10,12,14]), 	
\end{equation}
where
\begin{equation}
(x)_{0,\lambda}=1,\quad (x)_{n,\lambda}=x(x-\lambda)(x-2\lambda)\cdots(x-(n-1)\lambda),\quad (n\ge 1).\label{7}
\end{equation}
Note that $\lim_{\lambda\rightarrow 0}e_{\lambda}^{x}(t)=e^{xt}$.\\
In [10], the degenerate Stirling numbers of the second kind are given by
\begin{equation}
(x)_{n,\lambda}=\sum_{k=0}^{n}{n \brace k}_{\lambda}(x)_{k},\quad (n\ge 0). \label{8}
\end{equation}
Note that $\displaystyle\lim_{\lambda\rightarrow 0}{n\brace k}_{\lambda}={n \brace k}$. \par
Recently, the degenerate Whitney numbers of the second kind  $W_{m,\lambda}(n,k)$ are defined by
\begin{equation}
(mx+1)_{n,\lambda}=\sum_{k=0}^{n}W_{m,\lambda}(n,k)m^{k}(x)_{k},\quad (n\ge 0),\quad (\mathrm{see}\ [13]). \label{9}
\end{equation}
Note that $\lim_{\lambda\rightarrow 0}W_{m,\lambda}(n,k)=W_{m}(n,k)$. \\
In view of \eqref{4}, the degenerate Dowling polynomials are defined by
\begin{equation}
D_{m,\lambda}(n,x)=\sum_{k=0}^{n}W_{m,\lambda}(n,k)x^{k},\quad (n\ge 0),\quad (\mathrm{see}\ [13]). \label{10}
\end{equation}
From \eqref{9} and \eqref{10}, we can derive the following equations.
\begin{equation}
e_{\lambda}(t)\frac{1}{k!}\bigg(\frac{e_{\lambda}^{m}(t)-1}{m}\bigg)^{k}=\sum_{n=k}^{\infty}W_{m,\lambda}(n,k)\frac{t^{n}}{n!}, \label{11}	
\end{equation}
and
\begin{equation}
e_{\lambda}(t)e^{x\big(\frac{e_{\lambda}^{m}(t)-1}{m}\big)}=\sum_{n=0}^{\infty}D_{m,\lambda}(n,x)\frac{t^{n}}{n!},\quad (\mathrm{see}\ [13]). \label{12}
\end{equation}

For $r\ge 0$, the degenerate $r$-Whitney numbers of the second kind is given by
\begin{equation}
(mx+r)_{n,\lambda}=\sum_{k=0}^{n}W_{m,\lambda}^{(r)}(n,k)m^{k}(x)_{k},\quad (n\ge 0),\quad (\mathrm{see}\ [13]). \label{13}
\end{equation}
From \eqref{13}, we note that
\begin{equation}
e_{\lambda}^{r}(t)\frac{1}{k!}\bigg(\frac{e_{\lambda}^{m}(t)-1}{m}\bigg)^{k}=\sum_{n=k}^{\infty}W_{m,\lambda}^{(r)}(n,k)\frac{t^{n}}{n!},\quad (\mathrm{see}\ [13]). \label{14}
\end{equation}
The degenerate $r$-Dowling polynomials are defined by
\begin{equation}
e_{\lambda}^{r}(t)e^{x(\frac{e_{\lambda}^{m}(t)-1}{m})}=\sum_{n=0}^{\infty}D_{m,\lambda}^{(r)}(n,x)\frac{t^{n}}{n!},\quad (\mathrm{see}\ [15]). \label{15}
\end{equation}
By \eqref{14} and \eqref{15}, we get
\begin{equation}
D_{m,\lambda}^{(r)}(n,x)=\sum_{k=0}^{n}W_{m,\lambda}^{(r)}(n,k)x^{k},\quad (n\ge  0),\quad (\mathrm{see}\ [11,15]).\label{16}	
\end{equation} \par
For any integer $k\ge 0$, the partial Bell polynomials are given by
\begin{equation}
\frac{1}{k!}\bigg(\sum_{m=1}^{\infty}x_{m}\frac{t^{m}}{m!}\bigg)^{k}=\sum_{n=k}^{\infty}B_{n,k}(x_{1},x_{2},\dots,x_{n-k+1})\frac{t^{n}}{n!},\quad (\mathrm{see}\ [7,12]), \label{17}	
\end{equation}
\begin{equation}
\begin{aligned}
&B_{n,k}(x_{1},x_{2},\dots,x_{n-k+1}) \\
&=\sum_{\substack{l_{1}+l_{2}+\cdots+l_{n-k+1}=k\\ l_{1}+2l_{2}+\cdots+(n-k+1)l_{n-k+1}=n}}\frac{n!}{l_{1}!l_{2}!\cdots l_{n-k+1}!}\bigg(\frac{x_{1}}{1!}\bigg)^{l_{1}}\bigg(\frac{x_{2}}{2!}\bigg)^{l_{2}}\cdots \bigg(\frac{x_{n-k+1}}{(n-k+1)!}\bigg)^{l_{n-k+1}}.\label{18}
\end{aligned}	
\end{equation}
The complete Bell polynomials are defined by
\begin{equation}
\exp\bigg(\sum_{i=1}^{\infty}x_{i}\frac{t^{i}}{i!}\bigg)=\sum_{n=0}^{\infty}B_{n}(x_{1},x_{2},\dots,x_{n})\frac{t^{n}}{n!},\quad (\mathrm{see}\ [12]).\label{19}
\end{equation}
Thus, by \eqref{18} and \eqref{19}, we get
\begin{equation}
B_{n}(x_{1},x_{2},\dots,x_{n})=\sum_{k=0}^{n}B_{n,k}(x_{1},x_{2},\dots,x_{n-k+1}),\quad (n\ge 0).\label{20}
\end{equation} \par
Recently, for any random variable $Y$ satisfying \eqref{0}, the probabilistic degenerate Stirling numbers of the second kind associated with $Y$ are introduced by
\begin{equation}
\frac{1}{k!}\big(E[e_{\lambda}^{Y}(t)]-1\big)^{k}=\sum_{n=k}^{\infty}{n \brace k}_{Y,\lambda}\frac{t^{n}}{n!},\quad (k\ge 0),\quad (\mathrm{see}\ [16]).\label{23}
\end{equation}
Note that ${n \brace k}_{Y,\lambda}={n \brace k}_{\lambda}$, when $Y=1$. \par

\section{Probabilistic degenerate Dowling polynomials associated with random variables}
{\it{Throughout this section, we assume that $Y$ is a random variable such that the moment generating function of $Y$
\begin{equation*}
E[e^{tY}]=\sum_{n=0}^{\infty}E[Y^{n}]\frac{t^{n}}{n!},\,\,(|t|<r) \,\,\,\mathrm{exists\,\, for\,\, some} \,\,r>0, \quad (\mathrm{see}\ [1,3,4,16,21,24,25,26]).
\end{equation*}
We let $(Y_{j})_{j\ge 1}$ be a sequence of mutually independent copies of the random variable $Y$, and let
\begin{displaymath}
	S_{0}=0,\quad S_{k}=Y_{1}+Y_{2}+\cdots+Y_{k},\quad (k\in\mathbb{N}),\quad(\mathrm{see}\ [1,3,4,16,25,26]).
\end{displaymath}}} \par
In view of \eqref{11}, we define the probabilistic degenerate Whitney numbers of the second kind associated with $Y$ by
\begin{equation}
\frac{1}{k!}\bigg(\frac{E[e_{\lambda}^{mY}(t)]-1}{m}\bigg)^{k}e_{\lambda}(t)=\sum_{n=k}^{\infty}W_{m,\lambda}^{Y}(n,k)\frac{t^{n}}{n!},\quad (k\ge 0).\label{24}
\end{equation}
When $Y=1$, we have $W_{m,\lambda}^{Y}(n,k)=W_{m,\lambda}(n,k)$, \,\, (see \eqref{11}). \\
From \eqref{24}, we note taht
\begin{align}
\sum_{n=k}^{\infty}W_{m,\lambda}^{Y}(n,k)\frac{t^{n}}{n!}&=\frac{1}{m^{k}k!}\sum_{j=0}^{k}\binom{k}{j}(-1)^{k-j}\Big(E\big[e_{\lambda}^{mY}(t)\big]\Big)^{j}e_{\lambda}(t)\label{25} \\
&=\frac{1}{m^{k}k!}\sum_{j=0}^{k}\binom{k}{j}(-1)^{k-j}E\big[e_{\lambda}^{mY_{1}}(t)\big]\cdots E\big[e_{\lambda}^{mY_{j}}(t)\big]e_{\lambda}(t)\nonumber \\
&=\frac{1}{m^{k}k!}\sum_{j=0}^{k}\binom{k}{j}(-1)^{k-j}E\Big[e_{\lambda}^{m(Y_{1}+\cdots+Y_{j})+1}(t)\Big] \nonumber \\
&=\sum_{n=0}^{\infty}\frac{1}{m^{k}k!}\sum_{j=0}^{k}\binom{k}{j}(-1)^{k-j}E\Big[(mS_{j}+1)_{n,\lambda}\Big]\frac{t^{n}}{n!}. \nonumber
\end{align}
Therefore, by comparing the coefficients on both sides of \eqref{25}, we obtain the following theorem.
\begin{theorem}
For $n,k$ with $n\ge k\ge 0$, we have
\begin{displaymath}
W_{m,\lambda}^{Y}(n,k)=\frac{1}{m^{k}k!}\sum_{j=0}^{k}\binom{k}{j}(-1)^{k-j}E\big[(mS_{j}+1)_{m,\lambda}\big].
\end{displaymath}
\end{theorem}
Now, we define the operator $\triangle_{my}$ by
\begin{equation}
\begin{aligned}
&\triangle_{my}f(x)=f(x+my)-f(x),\\
& \triangle_{my_{1},my_{2},\dots,my_{k}}f(x)= \triangle_{my_{1}}\circ \triangle_{my_{2}}\circ\cdots\circ \triangle_{my_{k}}f(x),\quad (\mathrm{see}\ [4-16]).
\end{aligned}\label{26}
\end{equation}
From \eqref{26}, we can derive the following equation
\begin{equation}
\triangle_{my_{1},my_{2},\dots,my_{k}}e_{\lambda}^{x}(t)=\big(e_{\lambda}^{my_{k}}-1\big)\big(e_{\lambda}^{my_{k-1}}-1\big)\cdots\big(e_{\lambda}^{my_{1}}-1\big)e_{\lambda}^{x}(t). \label{27}
\end{equation}
Thus, by \eqref{24} and \eqref{27}, we get
\begin{align}
\sum_{n=k}^{\infty}W_{m,\lambda}^{Y}(n,k)\frac{t^{n}}{n!}&=\frac{1}{k!m^{k}}e_{\lambda}(t)\Big(E[e_{\lambda}^{mY}(t)]-1\Big)^{k} \label{28} \\
&=\frac{1}{k!m^{k}}e_{\lambda}(t)
\Big(E\big[e_{\lambda}^{mY_{1}}(t)\big]-1\Big) \Big(E\big[e_{\lambda}^{mY_{2}}(t)\big]-1\Big)\cdots \Big(E\big[e_{\lambda}^{mY_{k}}(t)\big]-1\Big)\nonumber \\
&=\frac{1}{k!m^{k}} E[\triangle_{mY_{1},mY_{2},\dots,Y_{k}}e_{\lambda}^{x}(t)]\vert_{x=1}=\sum_{n=0}^{\infty}\frac{1}{k!m^{k}}E[\triangle_{mY_{1},mY_{2},\dots,mY_{k}}(x)_{n,\lambda}] \vert_{x=1}\frac{t^{n}}{n!}.\nonumber	
\end{align}
By comparing the coefficients on both sides of \eqref{28}, we obtain the following theorem.
\begin{theorem}
For $n\ge k\ge 0$, we have
\begin{displaymath}
\frac{1}{k!m^{k}} E[\triangle_{mY_{1},mY_{2},\dots,mY_{k}}(1)_{n,\lambda}]=W_{m,\lambda}^{Y}(n,k),
\end{displaymath}
where $E[\triangle_{mY_{1},mY_{2},\dots,mY_{k}}(1)_{n,\lambda}]=E[\triangle_{mY_{1},mY_{2},\dots,mY_{k}}(x)_{n,\lambda}] \vert_{x=1}$, and we use the same convention as this one in the sequel.
\end{theorem}
From the equation (45) of [16], we note that
\begin{equation}
E\Big[\triangle_{mY_{1},mY_{2},\dots,mY_{k}}f(x)\Big]=\sum_{l=0}^{k}\binom{k}{l}(-1)^{k-l}E\Big[f(x+mS_{l})\Big]. \label{29}	
\end{equation}
Taking $f(x)=(x)_{n,\lambda}$ in \eqref{29} and then evaluating at $x=1$, we get
\begin{equation}
E\Big[\triangle_{mY_{1},mY_{2},\dots,mY_{k}}(1)_{n,\lambda}\Big]=\sum_{l=0}^{k}\binom{k}{l}(-1)^{k-l}E\Big[(mS_{l}+1)_{n,\lambda}\Big].\label{30}	
\end{equation}
Thus, by using \eqref{8}, \eqref{30} and Theorem 2.2, we get
\begin{align}
W_{m,\lambda}^{Y}(n,k)&=\frac{1}{k!m^{k}} E[\triangle_{mY_{1},mY_{2},\dots,mY_{k}}(1)_{n,\lambda}]=\frac{1}{k!m^{k}}\sum_{l=0}^{k}\binom{k}{l}(-1)^{k-l}E\Big[(mS_{l}+1)_{n,\lambda}\Big] \label{31} \\
&=\frac{1}{k!m^{k}}\sum_{l=0}^{k}\binom{k}{l}(-1)^{k-l}\sum_{j=0}^{n}{n \brace j}_{\lambda}E[(mS_{l}+1)_{j}] \nonumber \\
&=\frac{1}{k!m^{k}}\sum_{j=0}^{n}{n \brace j}_{\lambda}\sum_{l=0}^{k}\binom{k}{l}(-1)^{k-l}E\Big[(mS_{l}+1)_{j}\Big].\nonumber	
\end{align}
Therefore, by \eqref{31}, we obtain the following theorem.
\begin{theorem}
For $n\ge k\ge 0$, we have
\begin{displaymath}
W_{m,\lambda}^{Y}(n,k)=\frac{1}{m^{k}k!}\sum_{j=0}^{n}{n \brace j}_{\lambda}\sum_{l=0}^{k}\binom{k}{l}(-1)^{k-l}E\Big[(mS_{l}+1)_{j}\Big].
\end{displaymath}
\end{theorem}
Before proceeding further, we recall the well-known binomial inversion given by
\begin{equation}
a_{k}=\sum_{l=0}^{k}\binom{k}{l}b_{l} \Longleftrightarrow b_{k}=\sum_{l=0}^{k}(-1)^{k-l}\binom{k}{l}a_{l}. \label{31-1}
\end{equation}
By \eqref{30}, \eqref{31-1} and Theorem 2.2, we get
\begin{align}
&\sum_{k=0}^{N}E\Big[(mS_{k}+1)_{n,\lambda}	\Big]=\sum_{k=0}^{N}\sum_{l=0}^{k}\binom{k}{l}E\Big[\triangle_{mY_{1},mY_{2},\dots,mY_{l}}(1)_{n,\lambda}\Big] \label{32} \\
&\quad =\sum_{l=0}^{N}E\Big[\triangle_{mY_{1},mY_{2},\dots,mY_{l}}(1)_{n,\lambda}\Big]\sum_{k=l}^{N}\binom{k}{l}=\sum_{l=0}^{N}\binom{N+1}{l+1}E\Big[\triangle_{mY_{1},mY_{2},\dots,mY_{l}}(1)_{n,\lambda}\Big] \nonumber \\
&\quad =\sum_{l=0}^{N}l!m^{l}W_{m,\lambda}^{Y}(n,l)\binom{N+1}{l+1}. \nonumber
\end{align}
Therefore, by \eqref{32}, we obtain the following theorem.
\begin{theorem}
For $N\ge 0$, we have
\begin{displaymath}
\sum_{k=0}^{N}E\Big[(mS_{k}+1)_{n,\lambda}\Big]=\sum_{l=0}^{N}l!m^{l}\binom{N+1}{l+1}W_{m,\lambda}^{Y}(n,l).
\end{displaymath}
\end{theorem}
Now, we define the probabilistic degenerate Dowling polynomials associated with $Y$ by
\begin{equation}
D_{m,\lambda}^{Y}(n,x)=\sum_{k=0}^{n}W_{m,\lambda}^{Y}(n,k)x^{k},\quad (n\ge 0). \label{35}
\end{equation}
Note that if $Y=1$, then $D_{m,\lambda}^{Y}(n,x)=D_{m,\lambda}(n,x),\ (n\ge 0)$. When $x=1$, $D_{m,\lambda}^{Y}(n)=D_{m,\lambda}^{Y}(n,1)$ are called the probabilistic degenerate Dowling numbers associated with $Y$. \\
From \eqref{24} and \eqref{35}, we note that
\begin{align}
\sum_{n=0}^{\infty}D_{m,\lambda}^{Y}(n,x)\frac{t^{n}}{n!}&=\sum_{n=0}^{\infty}\sum_{k=0}^{n}W_{m,\lambda}^{Y}(n,k)x^{k}\frac{t^{n}}{n!}\label{36}\\
&=\sum_{k=0}^{\infty}x^{k}\sum_{n=k}^{\infty}W_{m,\lambda}^{Y}(n,k)\frac{t^{n}}{n!}=\sum_{k=0}^{\infty}x^{k}e_{\lambda}(t)\frac{1}{k!}\bigg(\frac{E[e_{\lambda}^{mY}(t)]-1}{m}\bigg)^{k}	\nonumber \\
&=e_{\lambda}(t)e^{x\big(\frac{E[e_{\lambda}^{mY}(t)]-1}{m}\big)}. \nonumber
\end{align}
Therefore, by \eqref{36}, we obtain the following theorem.
\begin{theorem}
Let $Y$ be a random variable. Then we have
\begin{equation}
e_{\lambda}(t)e^{x\big(\frac{E[e_{\lambda}^{mY}(t)]-1}{m}\big)}=\sum_{n=0}^{\infty}D_{m,\lambda}^{Y}(n,x)\frac{t^{n}}{n!}. \label{37}
\end{equation}
\end{theorem}
By \eqref{37}, we get
\begin{align}
&\sum_{n=0}^{\infty}D_{m,\lambda}^{Y}(n,x)\frac{t^{n}}{n!}=e_{\lambda}(t)e^{x\big(\frac{E[e_{\lambda}^{mY}(t)]-1}{m}\big)}=e^{-\frac{x}{m}}e_{\lambda}(t)e^{\frac{x}{m}E[e_{\lambda}^{mY}(t)]} \label{38}\\
&=e^{-\frac{x}{m}}e_{\lambda}(t)\sum_{k=0}^{\infty}\frac{\big(\frac{x}{m}\big)^{k}}{k!}\Big(E\big[e_{\lambda}^{mY}(t)\big]\Big)^{k}=e^{-\frac{x}{m}}\sum_{k=0}^{\infty}\frac{\big(\frac{x}{m}\big)^{k}}{k!}E\Big[e_{\lambda}^{m(Y_{1}+\cdots+Y_{k})+1}(t)\Big] \nonumber \\
&=\sum_{n=0}^{\infty}e^{-\frac{x}{m}}\sum_{k=0}^{\infty}\frac{x^{k}}{m^{k}k!}E\Big[(mS_{k}+1)_{n,\lambda}\Big]\frac{t^{n}}{n!}.\nonumber
\end{align}
Therefore, by \eqref{38}, we obtain the following Dobinski-like formula.
\begin{theorem}
For $n\ge 0$, we have
\begin{equation}
D_{m,\lambda}^{Y}(n,x)=e^{-\frac{x}{m}}\sum_{k=0}^{\infty}\frac{x^{k}}{m^{k}k!}E\Big[(mS_{k}+1)_{n,\lambda}\Big].\label{39}	
\end{equation}
\end{theorem}
From \eqref{37}, we note that
\begin{align}
&\sum_{n=0}^{\infty}D_{m,\lambda}^{Y}(n,x)\frac{t^{n}}{n!}=e_{\lambda}(t)e^{x\big(\frac{E[e_{\lambda}^{mY}(t)]-1}{m}\big)}=e_{\lambda}(t)e^{\frac{x}{m}\sum_{j=1}^{\infty}E[(mY)_{j,\lambda}]\frac{t^{j}}{j!}}\label{40} \\
&=e_{\lambda}(t)\sum_{k=0}^{\infty}\frac{1}{k!}\bigg(\frac{x}{m}\sum_{j=1}^{\infty}E[(mY)_{j,\lambda}]\frac{t^{j}}{j!}\bigg)^{k} \nonumber \\
&=e_{\lambda}(t)\sum_{k=0}^{\infty}\sum_{n=k}^{\infty}B_{n,k}\bigg(\frac{x}{m}E[(mY)_{1,\lambda}], \frac{x}{m}E[(mY)_{2,\lambda}],\dots, \frac{x}{m}E[(mY)_{n-k+1,\lambda}]\bigg)\frac{t^{n}}{n!}\nonumber \\
&=\sum_{j=0}^{\infty}\frac{(1)_{j,\lambda}}{j!}t^{j} \sum_{l=0}^{\infty}\sum_{k=0}^{l}B_{l,k} \bigg(\frac{x}{m}E[(mY)_{1,\lambda}], \frac{x}{m}E[(mY)_{2,\lambda}],\dots, \frac{x}{m}E[(mY)_{l-k+1,\lambda}]\bigg)\frac{t^{l}}{l!}\nonumber \\
&=\sum_{n=0}^{\infty}\sum_{l=0}^{n}\binom{n}{l}(1)_{n-l,\lambda}\sum_{k=0}^{l}B_{l,k} \bigg(\frac{x}{m}E[(mY)_{1,\lambda}], \frac{x}{m}E[(mY)_{2,\lambda}],\dots, \frac{x}{m}E[(mY)_{l-k+1,\lambda}]\bigg)\frac{t^{n}}{n!}.\nonumber
\end{align}
Therefore, by comparing the coefficients on both sides of \eqref{40}, we obtain the following theorem.
\begin{theorem}
For $n\ge 0$, we have
\begin{equation*}
\begin{aligned}
&D_{m,\lambda}^{Y}(n,x)\\
&\quad =\sum_{l=0}^{n}\sum_{n=0}^{l}\binom{n}{l}(1)_{n-l,\lambda}B_{l,k} \bigg(\frac{x}{m}E[(mY)_{1,\lambda}], \frac{x}{m}E[(mY)_{2,\lambda}],\dots, \frac{x}{m}E[(mY)_{l-k+1,\lambda}]\bigg).
\end{aligned}
\end{equation*}
\end{theorem}
From \eqref{37}, we have
\begin{align}
&\sum_{n=0}^{\infty}D_{m,\lambda}^{Y}(n+1,x)\frac{t^{n}}{n!}=\frac{d}{dt}\sum_{n=0}^{\infty}D_{m,\lambda}^{Y}(n,x)\frac{t^{n}}{n!}=\frac{d}{dt}\bigg(e_{\lambda}(t)e^{x\big(\frac{E[e_{\lambda}^{mY}(t)]-1}{m}\big)}\bigg) \label{41}\\
&=e_{\lambda}^{1-\lambda}(t) e^{x\big(\frac{E[e_{\lambda}^{mY}(t)]-1}{m}\big)}+\frac{xE[mYe_{\lambda}^{mY-\lambda}(t)]}{m}e_{\lambda}(t) e^{x\big(\frac{E[e_{\lambda}^{mY}(t)]-1}{m}\big)} \nonumber \\
&=\sum_{l=0}^{\infty}(-1)^{l}\lambda^{l}t^{l}\sum_{k=0}^{\infty}D_{m,\lambda}^{Y}(k,x)\frac{t^{k}}{k!}+\sum_{l=0}^{\infty}\frac{x}{m}E\big[(mY)_{l+1,\lambda}\big]\frac{t^{l}}{l!}\sum_{k=0}^{\infty}D_{m,\lambda}^{Y}(k,x)\frac{t^{k}}{k!}\nonumber\\
&=\sum_{n=0}^{\infty}\sum_{k=0}^{n}(-1)^{n-k}\lambda^{n-k}\frac{D_{m,\lambda}^{Y}(k,x)}{k!}t^{n}+\sum_{n=0}^{\infty}\sum_{k=0}^{n}\binom{n}{k}D_{m,\lambda}^{Y}(k,x)\frac{x}{m}E[(mY)_{n-k+1,\lambda}]\frac{t^{n}}{n!}\nonumber \\
&=\sum_{n=0}^{\infty}\bigg(\sum_{k=0}^{n}(-1)^{n-k}\lambda^{n-k}\frac{n!}{k!}D_{m,\lambda}^{Y}(k,x)+\frac{x}{m}\sum_{k=0}^{n}\binom{n}{k}D_{m,\lambda}^{Y}(k,x)E[(mY)_{n-k+1,\lambda}]\bigg)\frac{t^{n}}{n!}.\nonumber
\end{align}
Therefore, by comparing the coefficients on both sides of \eqref{41}, we obtain the following theorem.
\begin{theorem}
For $n\ge 0$, we have
\begin{equation*}
\begin{aligned}
&D_{m,\lambda}^{Y}(n+1,x)\\
&= \sum_{k=0}^{n}(-1)^{n-k}\lambda^{n-k}\frac{n!}{k!}D_{m,\lambda}^{Y}(k,x)+\frac{x}{m}\sum_{k=0}^{n}\binom{n}{k}D_{m,\lambda}^{Y}(k,x)E[(mY)_{n-k+1,\lambda}].
\end{aligned}
\end{equation*}
\end{theorem}
By \eqref{37}, we get
\begin{align}
&\sum_{n=0}^{\infty}\sum_{k=0}^{n}\binom{n}{k}(1)_{n-k,\lambda}D_{m,\lambda}^{Y}(k,x+y)\frac{t^{n}}{n!}=e_{\lambda}(t)e_{\lambda}(t) e^{(x+y)\big(\frac{E[e_{\lambda}^{mY}(t)]-1}{m}\big)}\label{42}\\
&=e_{\lambda}(t) e^{x\big(\frac{E[e_{\lambda}^{mY}(t)]-1}{m}\big)} e_{\lambda}(t) e^{y\big(\frac{E[e_{\lambda}^{mY}(t)]-1}{m}\big)}	\nonumber \\
&=\sum_{l=0}^{\infty}D_{m,\lambda}^{Y}(l,x)\frac{t^{l}}{l!}\sum_{k=0}^{\infty}D_{m,\lambda}^{Y}(k,y)\frac{t^{k}}{k!}=\sum_{n=0}^{\infty}\bigg(\sum_{k=0}^{n}\binom{n}{k}D_{m,\lambda}^{Y}(n-k,x)D_{m,\lambda}^{Y}(k,y)\bigg)\frac{t^{n}}{n!}.\nonumber
\end{align}
By \eqref{42}, we obtain the following theorem which shows that $D_{m,\lambda}^{Y}(n,x)$ does not satisfy the binomial identity.
\begin{theorem}
For $n\ge 0$, we have
\begin{displaymath}
\sum_{k=0}^{n}\binom{n}{k}(1)_{n-k,\lambda}D_{m,\lambda}^{Y}(k,x+y)=\sum_{k=0}^{n}\binom{n}{k}D_{m,\lambda}^{Y}(n-k,x)D_{m,\lambda}^{Y}(k,y).
\end{displaymath}
\end{theorem}
Let $r$ be a nonnegative integer. Then we consider the probabilistic degenerate $r$-Whitney numbers of the second kind associated with $Y$ defined by
\begin{equation}
\frac{1}{k!}\bigg(\frac{E[e_{\lambda}^{mY}(t)]-1}{m}\bigg)^{k}e_{\lambda}^{r}(t)=\sum_{n=k}^{\infty}W_{m,\lambda}^{(Y,r)}(n,k)\frac{t^{n}}{n!},\quad (k \ge 0).\label{43}
\end{equation}
When $Y=1,$ we have $W_{m,\lambda}^{(Y,r)}(n,k)=W_{m,\lambda}^{(r)}(n,k)$,\,\,(see \eqref{14}). \\
From \eqref{43}, we have
\begin{align}
\sum_{n=k}^{\infty}W_{m,\lambda}^{(Y,r)}(n,k)\frac{t^{n}}{n!}&=\frac{1}{k!m^{k}}\sum_{j=0}^{k}\binom{k}{j}(-1)^{k-j}\Big(E[e_{\lambda}^{mY}(t)]\Big)^{j}e_{\lambda}^{r}(t)\label{44}\\	
&=\frac{1}{m^{k}k!}\sum_{j=0}^{k}\binom{k}{j}(-1)^{k-j}E\Big[e_{\lambda}^{mS_{j}+r}(t)\Big]\nonumber \\
&=\sum_{n=0}^{\infty}\bigg(\frac{1}{m^{k}k!}\sum_{j=0}^{k}\binom{k}{j}(-1)^{k-j}E\Big[(mS_{j}+r)_{n,\lambda}\Big]\bigg)\frac{t^{n}}{n!}.\nonumber
\end{align}
Therefore, by comparing the coefficients on both sides of \eqref{44}, we obtain the following theorem.
\begin{theorem}
For $r\in\mathbb{N}\cup\{0\}$ and $n\ge k\ge 0$, we have
\begin{displaymath}
W_{m,\lambda}^{(Y,r)}(n,k)=\frac{1}{m^{k}k!}\sum_{j=0}^{k}\binom{k}{j}(-1)^{k-j}E\Big[(mS_{j}+r)_{n,\lambda}\Big].
\end{displaymath}
\end{theorem}
From \eqref{37}, we note that
\begin{align}
&\sum_{n=0}^{\infty}\sum_{k=0}^{n}\binom{n}{k}(x-1)_{n-k,\lambda}D_{m,\lambda}^{Y}(k,x)\frac{t^{n}}{n!}=\bigg(e_{\lambda}(t)e^{\big(\frac{E[e_{\lambda}^{mY}(t)]-1}{m}\big)}\bigg)^{x} \label{45}\\
&=\bigg(e_{\lambda}(t)e^{\big(\frac{E[e_{\lambda}^{mY}(t)]-1}{m}\big)}-1+1\bigg)^{x}=\sum_{k=0}^{\infty}\binom{x}{k}\bigg(e_{\lambda}(t) e^{\big(\frac{E[e_{\lambda}^{mY}(t)]-1}{m}\big)}-1\bigg)^{k}\nonumber \\
&=\sum_{k=0}^{\infty}\binom{x}{k}k!\frac{1}{k!}\bigg(\sum_{j=1}^{\infty}D_{m,\lambda}^{Y}(j)\frac{t^{j}}{j!}\bigg)^{k} \nonumber \\
&=\sum_{k=0}^{\infty}\binom{x}{k}k!\sum_{n=k}^{\infty}B_{n,k}\Big(D_{m,\lambda}^{Y}(1), D_{m,\lambda}^{Y}(2),\dots, D_{m,\lambda}^{Y}(n-k+1)\Big)\frac{t^{n}}{n!}\nonumber \\
&=\sum_{n=0}^{\infty}\sum_{k=0}^{n}\binom{x}{k}k! B_{n,k}\Big(D_{m,\lambda}^{Y}(1), D_{m,\lambda}^{Y}(2),\dots, D_{m,\lambda}^{Y}(n-k+1)\Big)\frac{t^{n}}{n!}.\nonumber
\end{align}
Therefore, by comparing the coefficients on both sides of \eqref{45}, we obtain the following theorem.
\begin{theorem}
For $n\ge 0$, we have
\begin{displaymath}
\sum_{k=0}^{n}\binom{n}{k}(x-1)_{n-k,\lambda}D_{m,\lambda}^{Y}(k,x)= \sum_{k=0}^{n}\binom{x}{k}k! B_{n,k}\Big(D_{m,\lambda}^{Y}(1), D_{m,\lambda}^{Y}(2),\dots, D_{m,\lambda}^{Y}(n-k+1)\Big).
\end{displaymath}
\end{theorem}
Now, we observe that
\begin{equation}
t e_{\lambda}(t)e^{x\big(\frac{E[e_{\lambda}^{mY}(t)]-1}{m}\big)}=t\sum_{j=0}^{\infty}D_{m,\lambda}^{Y}(j,x)\frac{t^{j}}{j!}=\sum_{j=1}^{\infty}jD_{m,\lambda}^{Y}(j-1,x)\frac{t^{j}}{j!}. \label{46}
\end{equation}
Thus, by \eqref{43} and \eqref{46}, we get
\begin{align}
&\bigg(\sum_{j=1}^{\infty}jD_{m,\lambda}^{Y}(j-1,x)\frac{t^{j}}{j!}\bigg)^{k}=t^{k}\bigg(e^{x\big(\frac{E[e_{\lambda}^{mY}(t)]-1}{m}\big)} e_{\lambda}(t)\bigg)^{k}\label{47} \\
&=t^{k}\sum_{j=0}^{\infty}k^{j}x^{j}\frac{1}{j!}\Big(\frac{E[e_{\lambda}^{mY}(t)]-1}{m}\Big)^{j}e_{\lambda}^{k}(t)=t^{k}\sum_{j=0}^{\infty}k^{j}x^{j}\sum_{n=j}^{\infty}W_{m,\lambda}^{(Y,k)}(n,j)\frac{t^{n}}{n!}\nonumber \\
&=\sum_{n=0}^{\infty}\sum_{j=0}^{n}k^{j}x^{j}W_{m,\lambda}^{(Y,k)}(n,j)\frac{t^{n+k}}{n!}=\sum_{n=k}^{\infty}k!\sum_{j=0}^{n-k}\binom{n}{k}k^{j}x^{j}W_{m,\lambda}^{(Y,k)}(n-k,j)\frac{t^{n}}{n!}. \nonumber
\end{align}
From \eqref{47}, we can derive the following equation
\begin{equation}
\begin{aligned}
&\sum_{n=k}^{\infty}\sum_{j=0}^{n-k}\binom{n}{k}k^{j}x^{j}W_{m,\lambda}^{(Y,k)}(n-k,j)\frac{t^{n}}{n!}=\frac{1}{k!}\bigg(\sum_{j=1}^{\infty}jD_{m,\lambda}^{Y}(j-1,x)\frac{t^{j}}{j!}\bigg)^{k}\\
&=\sum_{n=k}^{\infty}B_{n,k}\Big(D_{m,\lambda}^{Y}(0,x),2D_{m,\lambda}^{Y}(1,x),3D_{m,\lambda}^{Y}(2,x),\dots,(n-k+1)D_{m,\lambda}^{Y}(n-k,x)\Big)\frac{t^{n}}{n!}.
\end{aligned}\label{48}
\end{equation}
Therefore, by \eqref{48}, we obtain the following theorem.
\begin{theorem}
For $n,k$ with $n\ge k\ge 0$, we have
\begin{align*}
&\sum_{j=0}^{n-k}\binom{n}{k}k^{j}x^{j}W_{m,\lambda}^{(Y,k)}(n-k,j)\\
&= B_{n,k}\Big(D_{m,\lambda}^{Y}(0,x),2D_{m,\lambda}^{Y}(1,x),3D_{m,\lambda}^{Y}(2,x),\dots,(n-k+1)D_{m,\lambda}^{Y}(n-k,x)\Big).
\end{align*}
\end{theorem}
We note that
\begin{align}
&\sum_{n=k}^{\infty} B_{n,k}\Big(D_{m,\lambda}^{Y}(1,x)-(1)_{1,\lambda},D_{m,\lambda}^{Y}(2,x)-(1)_{2,\lambda},\dots,D_{m,\lambda}^{Y}(n-k+1,x)-(1)_{n-k+1,\lambda}\Big)\frac{t^{n}}{n!}\label{49} \\
&=\frac{1}{k!}\bigg(\sum_{j=1}^{\infty}\Big(D_{m,\lambda}^{Y}(j,x)-(1)_{j,\lambda}\Big)\frac{t^{j}}{j!}\bigg)^{k}=\frac{1}{k!}\bigg(e_{\lambda}(t)\label{49} e^{x\big(\frac{E[e_{\lambda}^{mY}(t)]-1}{m}\big)}-e_{\lambda}(t)\bigg)^{k}\nonumber \\
&=e_{\lambda}^{k}(t)\frac{1}{k!}\bigg(e^{x\big(\frac{E[e_{\lambda}^{mY}(t)]-1}{m}\big)}-1\bigg)^{k}=\sum_{j=k}^{\infty}{j \brace k}x^{j}\frac{1}{j!}\bigg(\frac{E[e_{\lambda}^{mY}(t)]-1}{m}\bigg)^{j}e_{\lambda}^{k}(t) \nonumber \\
&=\sum_{j=k}^{\infty}{j \brace k}x^{j}\sum_{n=j}^{\infty}W_{m,\lambda}^{(Y,k)}(n,j)\frac{t^{n}}{n!}=\sum_{n=k}^{\infty}\sum_{j=k}^{n}{j \brace k}W_{m,\lambda}^{(Y,k)}(n,j)x^{j}\frac{t^{n}}{n!}.\nonumber
\end{align}
Therefore, by \eqref{49}, we obtain the following theorem.
\begin{theorem}
For $n,k$ with $n\ge k\ge 0$, we have
\begin{equation*}
\begin{aligned}
&B_{n,k}\Big(D_{m,\lambda}^{Y}(1,x)-(1)_{1,\lambda},D_{m,\lambda}^{Y}(2,x)-(1)_{2,\lambda},\dots,D_{m,\lambda}^{Y}(n-k+1,x)-(1)_{n-k+1,\lambda}\Big)\\
&=\sum_{j=k}^{n}{j \brace k}W_{m,\lambda}^{(Y,k)}(n,j)x^{j}.
\end{aligned}
\end{equation*}
\end{theorem}
From \eqref{37}, we have
\begin{align}
&\sum_{n=0}^{\infty}\bigg(\frac{d}{dx}\bigg)^{k}D_{m,\lambda}^{Y}(n,x)\frac{t^{n}}{n!}=\bigg(\frac{d}{dx}\bigg)^{k}e^{x\big(\frac{E[e_{\lambda}^{mY}(t)]-1}{m}\big)}e_{\lambda}(t) \label{50}\\
&=k!\frac{1}{k!}\bigg(\frac{E[e_{\lambda}^{mY}(t)]-1}{m}\bigg)^{k} e^{x\big(\frac{E[e_{\lambda}^{mY}(t)]-1}{m}\big)}e_{\lambda}(t)\nonumber \\
&=\frac{k!}{m^{k}}\frac{1}{k!}\Big(E[e_{\frac{\lambda}{m}}^{Y}(mt)]-1\Big)^{k} e^{x\big(\frac{E[e_{\lambda}^{mY}(t)]-1}{m}\big)}e_{\lambda}(t) \nonumber \\
&=\frac{k!}{m^{k}}\sum_{l=k}^{\infty}{l \brace k}_{Y,\frac{\lambda}{m}}\frac{m^{l}t^{l}}{l!}\sum_{j=0}^{\infty}D_{m,\lambda}^{Y}(j,x)\frac{t^{j}}{j!} \nonumber \\
&=\sum_{n=k}^{\infty}k!\sum_{j=0}^{n-k}\binom{n}{j}D_{m,\lambda}^{Y}(j,x){n-j \brace k}_{Y,\frac{\lambda}{m}}m^{n-k-j}\frac{t^{n}}{n!}.\nonumber
\end{align} \\
In particular, for $k=1$, we have
\begin{align}
&\sum_{n=1}^{\infty}\frac{d}{dx}D_{m,\lambda}^{Y}(n,x)\frac{t^{n}}{n!}=\bigg(\frac{E[e_{\lambda}^{mY}(t)]-1}{m}\bigg) e^{x\big(\frac{E[e_{\lambda}^{mY}(t)]-1}{m}\big)} e_{\lambda}(t) \label{52}\\
&=\frac{1}{m}\sum_{l=1}^{\infty}E\big[(Y)_{l,\frac{\lambda}{m}}\big]\frac{m^{l}t^{l}}{l!} \sum_{j=0}^{\infty}D_{m,\lambda}^{Y}(j,x)\frac{t^{j}}{j!}\nonumber\\
&=\sum_{n=1}^{\infty}\sum_{j=0}^{n-1}\binom{n}{j}E[(Y)_{n-j,\frac{\lambda}{m}}]m^{n-j-1}D_{m,\lambda}^{Y}(j,x)\frac{t^{n}}{n!}. \nonumber
\end{align}\\
By comparing the coefficients on both sides of \eqref{50} and \eqref{52}, we obtain the following theorem.
\begin{theorem}
For $n,k\in\mathbb{N}$ with $n\ge k$, we have
\begin{displaymath}
\bigg(\frac{d}{dx}\bigg)^{k}D_{m,\lambda}^{Y}(n,x)=k!\sum_{j=0}^{n-k}\binom{n}{j}D_{m,\lambda}^{Y}(j,x){n-j \brace k}_{Y,\frac{\lambda}{m}}m^{n-k-j}.
\end{displaymath}
In particular, for $k=1$, we have
\begin{displaymath}
\frac{d}{dx}D_{m,\lambda}^{Y}(n,x)=\sum_{j=0}^{n-1}\binom{n}{j}E\big[(Y)_{n-j,\frac{\lambda}{m}}\big]D_{m,\lambda}^{Y}(j,x)m^{n-j-1},\quad (n \ge 1).
\end{displaymath}
\end{theorem}
We observe that
\begin{align}
&\sum_{n=k}^{\infty}W_{m,\lambda}^{Y}(n,k)\frac{t^{n}}{n!}=\frac{1}{k!}\bigg(\frac{E[e_{\lambda}^{mY}(t)]-1}{m}\bigg)^{k}e_{\lambda}(t) \label{54} \\
&=\frac{1}{m^{k}}\frac{1}{k!}\Big(E[e_{\frac{\lambda}{m}}^{Y}(mt)]-1\Big)^{k}e_{\lambda}(t)=\frac{1}{m^{k}}\frac{1}{k!}\bigg(\sum_{j=1}^{\infty}E[(Y)_{j,\frac{\lambda}{m}}]m^{j}\frac{t^{j}}{j!}\bigg)^{k}e_{\lambda}(t)\nonumber \\
&=\frac{1}{m^{k}}\sum_{l=k}^{\infty}B_{l,k}\Big(E[(Y)_{1,\frac{\lambda}{m}}]m, E[(Y)_{2,\frac{\lambda}{m}}]m^{2},\dots, E[(Y)_{l-k+1,\frac{\lambda}{m}}]m^{l-k+1}\Big)\frac{t^{l}}{l!} \sum_{j=0}^{\infty} (1)_{j,\lambda}\frac{t^{j}}{j!}\nonumber \\
&=\sum_{n=k}^{\infty}\frac{1}{m^k}\sum_{l=k}^{n}\binom{n}{l}B_{l,k}\Big(E[(Y)_{1,\frac{\lambda}{m}}]m, E[(Y)_{2,\frac{\lambda}{m}}]m^{2},\dots, E[(Y)_{l-k+1,\frac{\lambda}{m}}]m^{l-k+1}\Big)(1)_{n-l,\lambda}\frac{t^n}{n!}. \nonumber
\end{align}
Thus, by \eqref{54}, we get the next theorem.
\begin{theorem}
For $n,k$ with $n\ge k\ge 0$, we have
\begin{displaymath}
W_{m,\lambda}^{Y}(n,k)= \frac{1}{m^k}\sum_{l=k}^{n}\binom{n}{l}B_{l,k}\Big(E[(Y)_{1,\frac{\lambda}{m}}]m, E[(Y)_{2,\frac{\lambda}{m}}]m^{2},\dots, E[(Y)_{l-k+1,\frac{\lambda}{m}}]m^{l-k+1}\Big)(1)_{n-l,\lambda}.
\end{displaymath}
\end{theorem}
Now, we define the probabilistic degenerate $r$-Dowling polynomials associated with $Y$ by
\begin{equation}
D_{m,\lambda}^{(Y,r)}(n,x)=\sum_{k=0}^{n}W_{m,\lambda}^{(Y,r)}(n,k)x^{k},\quad (n\ge 0).\label{56}
\end{equation}
From \eqref{56}, we note that
\begin{align}
\sum_{n=0}^{\infty}D_{m,\lambda}^{(Y,r)}(n,x)\frac{t^{n}}{n!}&=\sum_{n=0}^{\infty}\sum_{k=0}^{n}W_{m,\lambda}^{(Y,r)}(n,k)x^{k}\frac{t^{n}}{n!} \label{57} \\
&=\sum_{k=0}^{\infty}x^{k}\sum_{n=k}^{\infty}W_{m,\lambda}^{(Y,r)}(n,k)\frac{t^{n}}{n!}=\sum_{k=0}^{\infty}x^{k}e_{\lambda}^{r}(t)\frac{1}{k!}\bigg(\frac{E[e_{\lambda}^{mY}(t)]-1}{m}\bigg)^{k}\nonumber \\
&=e_{\lambda}^{r}(t)	e^{x\big(\frac{E[e_{\lambda}^{mY}(t)]-1}{m}\big)}.\nonumber
\end{align}
By \eqref{57}, we obtain the following theorem.
\begin{theorem}
Let $r$ be a nonnegative integer. Then we have
\begin{equation}
e_{\lambda}^{r}(t) e^{x\big(\frac{E[e_{\lambda}^{mY}(t)]-1}{m}\big)}=\sum_{n=0}^{\infty}D_{m,\lambda}^{(Y,r)}(n,x)\frac{t^{n}}{n!}.\label{58}
\end{equation}
\end{theorem}
By \eqref{58} , we get
\begin{align}
\sum_{n=0}^{\infty}D_{m,\lambda}^{(Y,r)}(n,x)\frac{t^{n}}{n!}&=e_{\lambda}^{r}(t) e^{x\big(\frac{E[e_{\lambda}^{mY}(t)]-1}{m}\big)}\label{59} \\
&=e^{-\frac{x}{m}}e_{\lambda}^{r}(t)\sum_{k=0}^{\infty}\frac{\big(\frac{x}{m}\big)^{k}}{k!}\Big(E[e_{\lambda}^{mY}(t)]\Big)^{k}=e^{-\frac{x}{m}}\sum_{k=0}^{\infty}\frac{\big(\frac{x}{m}\big)^{k}}{k!}E[e_{\lambda}^{mS_{k}+r}]\nonumber \\
&=\sum_{n=0}^{\infty}e^{-\frac{x}{m}}\sum_{k=0}^{\infty}\frac{x^{k}}{m^{k}k!}E\big[(mS_{k}+r)_{n,\lambda}\big]\frac{t^{n}}{n!}.\nonumber
\end{align}
By \eqref{59}, we obtain the following Dobinski-like theorem.
\begin{theorem}
	For $n\ge 0$, we have
	\begin{displaymath}
		D_{m,\lambda}^{(Y,r)}(n,x)= e^{-\frac{x}{m}}\sum_{k=0}^{\infty}\frac{x^{k}}{k!m^{k}}E\Big[(mS_{k}+r)_{n,\lambda}\Big].
	\end{displaymath}
\end{theorem}
When $Y=1$, we have $D_{m,\lambda}^{(Y,r)}(n,x)=D_{m,\lambda}^{(r)}(n,x),\quad (n\ge 0).$

\section{Conclusion}
In this paper, we studied by using generating functions the probabilistic degenerate Whitney numbers of the second kind associated with $Y$ and the probabilistic degenerate Dowling polynomials associated with $Y$. Here $Y$ is a random variable satisfying the moment condition in \eqref{0}. In addition, we investigated generalizations of those numbers and polynomials, namely the probabilistic degenerate $r$-Whitney numbers of the second kind associated with $Y$ and the probabilistic degenerate $r$-Dowling polynomials associated with $Y$. We derived some properties, explicit expressions, certain identities, recurrence relations and generating functions for those numbers and polynomials. \par
It is one of our future projects to continue to study degenerate versions, $\lambda$-analogues and probabilistic versions of many special polynomials and numbers and to find their applications to physics, science and engineering as well as to mathematics.

\end{document}